\input form

\def\ee{{e}}
\def\realpart{\mathop{\rm Re}\nolimits}
\def\imaginary{\mathop{\rm Im}\nolimits}
\def\imunit{{\bf i}}
\def\der{{\rm d}}
\def\Lie{{\Lscr}}
\def\telchi{{T}}
\def\parder#1#2{{\partial#1 \over \partial#2}}
\def\scalprod#1#2{#1\cdot#2}
\def\zconj{\imunit\bar z}

\cita{BenGalGio-1989}
{G. Benettin, L. Galgani, A. Giorgilli: {\it Realization of holonomic
constraints and freezing of high frequency degrees of freedom in the
light of classical perturbation theory, part II}, Comm. Math. Phys,
{\bf 121}, 557--601 (1989).}
\cita{CecBig-2013b}
{M. Ceccaroni, J.D. Biggs: {\it Analytical perturbative method for
frozen orbits around the asteroid 433 Eros}, IAC2012 International
Astronautical Congress - Naples, IAC-12,C1,7,6,x14267 (2012).}
\cita{CecBig-2013a}
{M. Ceccaroni, J.D. Biggs: {\it Analytic perturbative theories in highly
inhomogeneous gravitational fields}, Icarus, {\bf 224}, 74--85 (2013).}
\cita{CecBisBig-2014}
{M. Ceccaroni, F. Biscani, J.D. Biggs: {\it Analytical method for
perturbed frozen orbit around an Asteroid in highly inhomogeneous
gravitational fields: a first approach}, Solar Syst. Res., {\bf
48}, 33--47 (2014).}
\cita{DepPalDep-2001}
{A. Deprit, J. Palaci\'an, E. Deprit: {\it The Relegation Algorithm},
 CeMDA, {\bf 79}, 157--182 (2001).}
\cita{FenNooVisYua-2015}
{J. Feng, R. Noomen, P.N. Visser, J. Yuan: {\it Modelling and analysis
 of periodic orbits around a contact binary asteroid},
 Astrophys. Space Sci., {\bf 357}, 1--18 (2015).}
\cita{Giorgilli-2003}
{A. Giorgilli: {\it Notes on exponential stability of Hamiltonian
systems}, in {\it Dynamical Systems, Part
I\/}, Pubbl. Cent. Ric. Mat. Ennio De Giorgi, Sc. Norm. Sup.  Pisa,
87--198 (2003).}
\cita{GioLocSan-2009}
{A. Giorgilli, U. Locatelli, M. Sansottera: {\it Kolmogorov and
Nekhoroshev theory for the problem of three bodies}, CeMDA, {\bf 104},
159--173 (2009).}
\cita{GioLocSan-2014}
{A. Giorgilli, U. Locatelli, M. Sansottera: {\it On the convergence of
an algorithm constructing the normal form for elliptic lower
dimensional tori in planetary systems}, CeMDA, {\bf 119}, 397--424
(2014).}
\cita{GioSan-2012}
{A. Giorgilli, M. Sansottera: {\it Methods of algebraic
manipulation in perturbation theory}, Workshop Series of the
Asociacion Argentina de Astronomia, {\bf 3}, 147--183 (2011).}
\cita{Grobner-60}
{W. Gr\"obner: {\it Die Lie-Reihen und Ihre Anwendungen}, Springer
Verlag, Berlin (1960); Italian transl.: {\it Le serie di Lie e le loro
applicazioni}, Cremonese, Roma (1973).}
\cita{LarSanLop-2013}
{M. Lara, J.F. San-Juan, L.M. L\'opez-Ochoa: {\it Averaging Tesseral
Effects: Closed Form Relegation versus Expansions of Elliptic Motion},
Math. Probl. Eng., {\bf 2013}, 570127 (2013).}
\cita{LarSanLop-2014}
{M. Lara, J.F. San-Juan, L.M. L\'opez-Ochoa: {\it Delaunay variables
approach to the elimination of the perigee in Artificial Satellite
Theory}, CeMDA, {\bf 120}, 39--56 (2014).}
\cita{Nekhoroshev-1977}
{N.N. Nekhoroshev: {\it Exponential estimates of the stability time of
near-integrable Hamiltonian systems}. English translation: Russ. Math.
Surveys, {\bf 32}, 1 (1977).}
\cita{Nekhoroshev-1979}
{N.N. Nekhoroshev: {\it Exponential estimates of the stability time of
near-integrable Hamiltonian systems, 2.} Trudy Sem.\ Im.\ G.\
Petrovskogo, {\bf 5}, 5 (1979).  English translation: {\it Topics in
modern Mathematics}, Petrovskij Semin., {\bf 5}, 1--58 (1985).}
\cita{NouTsiTzi-2015}
{A. Noullez, K. Tsiganis, S. Tzirti: {\it Satellite orbits design
using frequency analysis}, Adv. Space Res., {\bf 56},
163--175 (2015).}
\cita{Palacian-1992}
{J. Palaci\'an: {\it Teor{\'\i}a del Sat\'elite Artificial:
Arm\'onicos Teserales y su Relegaci\'on Mediante Simplificaciones
Algebraicas}, Ph.D. thesis, Universidad de Zaragoza (1992).}
\cita{Palacian-2002}
{J. Palaci\'an: {\it Normal Forms for Perturbed Keplerian Systems},
J. Differ. Equations, {\bf 180}, 471--519 (2002).}
\cita{ParMorKug-2014}
{P.C.P.M. Pardal, R. V. de Moraes, H. K. Kuga, {\it Effects of
Geopotential and Atmospheric Drag Effects on Frozen Orbits Using
Nonsingular Variables.}  Math. Probl. Eng. (2014).}
\cita{SanLhoLem-2014}
{M. Sansottera, C. Lhotka, A. Lema{\^\i}tre: {\it Effective stability
around the Cassini state in the spin-orbit problem}, CeMDA, {\bf 119},
75--89 (2014).}
\cita{SanLhoLem-2015}
{M. Sansottera, C. Lhotka, A. Lema{\^\i}tre: {\it Effective resonant
stability of Mercury}, MNRAS, {\bf 452}, 4145--4152 (2015).}
\cita{SanLocGio-2011}
{M. Sansottera, U. Locatelli, A. Giorgilli: {\it A Semi-Analytic
Algorithm for Constructing Lower Dimensional Elliptic Tori in
Planetary Systems}, CeMDA, {\bf 111}, 337--361 (2011).}
\cita{SanLocGio-2013}
{M. Sansottera, U. Locatelli, A. Giorgilli: {\it On the stability of
the secular evolution of the planar Sun-Jupiter-Saturn-Uranus system},
Math. Comput. Simul., {\bf 88}, 1--14 (2013).}

\def\testos{M. Sansottera and M. Ceccaroni}
\def\testod{Rigorous estimates for the relegation algorithm}

\title{Rigorous estimates for\riga the relegation algorithm}

\author{\it MARCO SANSOTTERA\hfil\break
Dipartimento di Matematica,
Universit\`a degli Studi di Milano,\hfil\break
via Saldini 50, 20133\ ---\ Milano, Italy.}

\author{\it MARTA CECCARONI\hfil\break
Dipartimento di Matematica,
Universit\`a degli Studi di Roma ``Tor Vergata'',\hfil\break
via della Ricerca Scientifica 1, 00133\ ---\ Roma, Italy.}

\abstract{%
We revisit the relegation algorithm by Deprit et al. (2001)
in the light of the rigorous Nekhoroshev's like theory.  This
relatively recent algorithm is nowadays widely used for implementing
closed form analytic perturbation theories, as it generalises the
classical Birkhoff normalisation algorithm.  The algorithm, here
briefly explained by means of Lie transformations, has been so far
introduced and used in a {\it formal} way, i.e. without providing any
rigorous convergence or asymptotic estimates. The overall aim of this
paper is to find such quantitative estimates and to show how the
results about stability over exponentially long times can be recovered
in a simple and effective way, at least in the non-resonant case.}

\section{sec:intro}{Introduction}
The relegation algorithm, firstly introduced by
Palaci\'an\bibref{Palacian-1992} and fully stated
in~\dbiref{DepPalDep-2001}, is an algorithm suited for canonical
simplifications of Hamiltonian problems and represents an extension of
the classical Birkhoff normal form method.  Given a system described
by a perturbed Hamiltonian, this procedure aims to {\it relegate the
effect of a desired part} of the perturbation to an arbitrarily small
reminder, by constructing a suitable sequence of canonical changes of
coordinates.  The computation of such a transformed Hamiltonian would
require an infinite sequence of canonical transformations which, in
general, cannot be done since the resulting series turns out to be
divergent.  However, the asymptotic character of the transformation
allows to progressively lower the influence of the so-called
remainder, i.e. the not yet relegated part of the Hamiltonian, by
considering a finite sequence of canonical transformations.  This is
very useful in practical applications, where one proceeds
order-by-order until the remainder can be considered, in some sense,
negligible.

Consider a Hamiltonian expanded in power series of a small parameter,
namely $H = H_0 + \epsilon H_1 + \ldots\ $, and in particular the case
$H_0=h_0+f_0$ where $h_0$ is an integrable Hamiltonian and $f_0$ is
for the moment a generic function.  The relegation aims at
constructing a first integral that at first order is essentially
$h_0$.  This is obtained by constructing a canonical change of
coordinates that puts the transformed Hamiltonian in a suitable {\it
normal form}. In this respect the relegation procedure is different
from that of Birkhoff, in which the normal form Hamiltonian is
requested to commute with $H_0$.  Indeed the effect of the function
$f_0$ is reduced to a small contribution, i.e., {\it relegated}.  The
relegation algorithm has been successfully used in celestial
mechanics\bibref{Palacian-2002}\bibref{CecBig-2013b}\bibref{CecBisBig-2014}
and artificial satellite
theory\bibref{LarSanLop-2013}\bibref{LarSanLop-2014}, with a
particular focus on the dynamics close to asteroids (the so-called
fast-rotating
case) \bibref{CecBig-2013a}\bibref{ParMorKug-2014}\bibref{FenNooVisYua-2015}\bibref{NouTsiTzi-2015}.
Let us remark that a similar normal form approach has also been
introduced in~\dbiref{BenGalGio-1989}, where the authors studied the
problem of the energy exchanges between a system of uncoupled harmonic
oscillators and a generic other dynamical system, playing the role of
the functions $h_0$ and $f_0$ of the relegation algorithm,
respectively.

Following the usual tradition in celestial mechanics, the relegation
algorithm has been introduced in a {\it formal} way.  However, to our
knowledge, the literature lacks of rigorous {\it convergence} or, at
least, asymptotic estimates for the algorithm.  We plan to do this in
the present paper.

\subsection{sbs:intro-relegation}{The relegation algorithm}
Throughout the paper we will make a wide use of the formalism of Lie
transforms which we briefly recall here.  We refer
to~\dbiref{Grobner-60} and~\dbiref{Giorgilli-2003} for an exhaustive
introduction.

The Lie transform $\telchi_{\Xscr}g$ of a generic function $g$ is defined as
$$
\telchi_{\Xscr}g=\sum_{j=0}^{\infty}E_{j}\,g
\qquad
{\rm with}
\qquad
E_{0}\,g=g\ ,
\qquad
E_{j}\,g=\sum_{i=1}^{j}\frac{i}{j}\Lie_{\Xscr_i}E_{j-i}\,g
\ ,
$$
where $\Xscr=\{\Xscr_1,\,\Xscr_2,\ldots\}$ is a sequence of generating
functions and $\lie{\Xscr}g$ is the Lie derivative of $g$ with
respect to $X$, i.e., the Poisson bracket $\{\Xscr,g\}$.  There is
also an explicit formula for the inverse, namely
$$
\telchi_{\Xscr}^{-1}g=\sum_{j=0}^{\infty}D_{j}\,g
\qquad
{\rm with}
\qquad
D_{0}\,g=g\ ,
\qquad
D_{j}\,g=-\sum_{i=1}^{j}\frac{i}{j}D_{j-i}\Lie_{\Xscr_i}g\ .
$$

The Lie transform is a generalisation of the Lie series, extremely
useful in perturbation theory as it can represent every near the
identity canonical transformation.  Moreover Lie transforms, being
defined in a recurrent explicit formula, are well suited to develop
effective algorithms and perform computations using computer algebra
(see, e.g., \dbiref{GioSan-2012}).

Consider a Hamiltonian
$$
H = \sum_{s\geq0} H_s\ ,
$$
where $H_s$ is a term of order $s$ in some small parameter.  As a
general fact in perturbation theory one aims at transforming the
Hamiltonian to a so-called normal form that we will denote as
$$
Z = \sum_{s\geq0} Z_s\ .
$$
The condition of being in normal form is that $Z$ should commute with
a certain function $h_0$, i.e. $\{Z,h_0\}=0$.  For instance, in
Birkhoff normalisation $Z$ must commute with $h_0=H_0$.  The algorithm
is developed by solving the equation
$$
\telchi_{\Xscr}Z=H\ ,
$$
for the
sequence of generating functions, $\Xscr$, and the normalised Hamiltonian $Z$.
In particular, at each order, we have to solve the so-called
homological equation
$$
\lie{H_0} \Xscr_s + Z_s = \Psi_s\ ,
$$
where $\Psi_s$ is a known function collecting all the terms of order
$s$, while $\Xscr_s$ and $Z_s$ are the unknowns to be determined.

In contrast with normalisation, in the relegation procedure $H_0$ is
split in two functions, i.e., $H_0 = h_0 + f_0$, also asking
$\{h_0,f_0\}=0$ and require that the sole $h_0$ is integrable.  Thus,
in the relegation algorithm, the normal form terms $Z_s$, with
$Z_0=h_0+f_0$, must commute with $h_0$ so that
it becomes a formal first integral for the normal form.  In this case, at each order,
the homological equation takes the form
$$
\lie{h_0} \Xscr_s + Z_s = \Psi_s\ ,
\formula{frm:hom2}
$$
where, again, $\Psi_s$ is a known function collecting all the terms of
order $s$, while $\Xscr_s$ and $Z_s$ are the unknowns to be
determined.
The crucial difference with respect to the Birkhoff
normalisation algorithm is that the generating sequence, $\Xscr$, is
determined using only the function $h_0$ and not the whole zero-order
term $H_0 = h_0+f_0$.  This tiny difference can have great impact in
specific problems arising in celestial mechanics and astrodynamics, as
illustrated in the papers quoted above.

As a matter of fact, in the relegation algorithm the homological
equation~\frmref{frm:hom2} is solved via a recurrence formula that aims
to counteract the effect of the terms generated by $f_0$.  Thus $f_0$
plays a special role, as detailed in section~\secref{sec:formale}.
Again, due to lack of convergence we can only {\it relegate} the
action of $f_0$ to a suitable order, thus making it as small as
possible.

\subsection{sbs:intro-risultati}{Statement of the results}
We collect here our main results, i.e., the asymptotic properties of
the relegation algorithm, together with estimates about the long time
stability of the approximate first integrals.  Precisely, we give
rigorous bounds for a truncated sequence of generating functions
$\Xscr$ and the corresponding transformed Hamiltonian $Z$.  All details
on the relegation algorithm and the quantitative estimates are included in
section~\secref{sec:formale} and~\secref{sec:stime}, respectively.

Consider a system of differential equations with Hamiltonian
$$
H(p,q,z,\zconj) = h_0(p) + \mu f_0(p,q,z,\zconj) + \epsilon H_1(p,q,z,\zconj)\ ,
\formula{frm:H0}
$$
with action-angle variables $p\in\Gscr \subseteq \reali^{n_1}$,
$q\in\toro^{n_1}$, and conjugate canonical variables
$(z,\zconj)\in\Bscr \subseteq \complessi^{2 n_2}$, where both $\Gscr$
and $\Bscr$ are open sets containing the origin and $n_1$, $n_2$ are
positive integers.  The quantities $\mu,\epsilon\in\reali$ are two
small parameters, with $\mu>\epsilon$.  The
latter request is not essential for the correctness of the proof.
However, in case $\mu<\epsilon$, the term $\mu f_0$ can be moved to
the perturbation and one can proceed by means of standard
normalisation procedure (see section~4 in~\dbiref{DepPalDep-2001}).

We consider the domain $\Dscr = \Gscr \times \toro^{n_1} \times \Bscr$
and introduce the extended domains\footnote{1}{Precisely,
$\Gscr_{\rho}=\big\{z\in\complessi^{n_1}:\max_{1\le j\le
n_1}|z_j|<\rho\big\}$,
$\toro^{n_1}_{\sigma}=\big\{q\in\complessi^{n_1}:\realpart
q_j\in\toro,\break\ \max_{1\le j\le n_1}|\imaginary
q_j|<\sigma\big\}\,$, $\Bscr_R=\{z\in\complessi^{2n_2}: \max_{1\le
j\le 2n_2} |z_j|<R\,\}$.} $\Dscr_{\rho,\sigma,R}
= \Gscr_{\rho} \times \toro^{n_1}_{\sigma} \times \Bscr_{R}\,$, where
$\Gscr_{\rho}\subset\complessi^{n_1}$ and
$\Bscr_R\subset\complessi^{2n_2}$ are complex open balls centred at
the origin with radii $\rho$ and $R$, respectively, while the
subscript $\sigma$, with $\sigma\in\reali$ such that $\sigma>0$,
denotes the usual complex extension of the torus.

Let us consider a generic analytic function
$g:\Dscr_{\rho,\sigma,R}\to\complessi$,
$$
g(p,q,z,\zconj) =
\sum_{{\scriptstyle{k\in\interi^{n_1}}}} g_{k}(p,z,\zconj)
e^{\imunit k\cdot q}\ ,
\formula{frm:funz}
$$
where $g_{k}:\Gscr_{\rho}\times \Bscr_{R}\to\complessi\,$.
We define the supremum norm
$$
\left|g\right|_{\rho,\sigma,R}=
\sup_{{\scriptstyle{p\in\Gscr_{\rho}}},\ {\scriptstyle{q\in\toro^{n_1}_{\sigma}}}
      \atop{{\scriptstyle{(z,\zconj)\in \Bscr_{R}}}}}
\big|g(p,q,z,\zconj)\big|\ .
$$
and the weighted Fourier norm
$$
\|g\|_{\rho,\sigma,R}=\sum_{{\scriptstyle{k\in\interi^{n_1}}}}
\big|g_{k}\big|_{\rho,R} e^{|k|\sigma}\ ,
\formula{frm:normafou}
$$
where
$$
\left|g_{k}\right|_{\rho,R}=
\sup_{{\scriptstyle{p\in\Gscr_{\rho}}}
      \atop{{\scriptstyle{(z,\zconj)\in \Bscr_{R}}}}}
\big|g_{k}(p,z,\zconj)\big|\ .
$$
Hereafter, we use the shorthand notations $|\cdot|_{\alpha}\,$ and
$\|\cdot\|_{\alpha}\,$ for $|\cdot|_{\alpha(\rho,R)}$ and
$\|\cdot\|_{\alpha(\rho,\sigma,R)}\,$.

\smallskip
The Hamiltonian~\frmref{frm:H0} is characterised by 8 real parameters,
namely $\rho$, $\sigma$, $R$ (analytic parameters), $\epsilon$, $\mu$
(perturbation parameters) and $\omega$, $\gamma$, $\tau$ (frequency
parameters), where the latter are introduced below.  We make the following
hypotheses:

\item{(i)}
$h_0(p)= \scalprod{\omega}{p}$, where $\omega\in\reali^{n_1}$ is a
fixed frequency vector;

\item{(ii)}
$h_0$, $f_0$ and $H_1$ are holomorphic bounded functions on the
extended domain $\Dscr_{\rho,2\sigma,R}$ with
$\|f_0\|_{\rho,2\sigma,R} \leq G$ for some positive real $G$;

\item{(iii)}
the functions $h_0$ and $f_0$ commute, i.e., $\poisson{h_0}{f_0}=0$.

Moreover, we introduce the so-called {\it resonance module} $\Mscr_\omega$ associated to the fixed frequency vector $\omega\in\reali^{n_1}$ as
$$
\Mscr_{\omega} = \{ k\in\interi^{n_1}: \scalprod{k}{\omega}=0\}\ .
$$
We assume the following additional hypothesis:

\item{(iv)}
for every positive integers $r$ and $K$, the frequency vector
$\omega\in\reali^{n_1}$ satisfies the Diophantine condition
$$
\alpha_r = \min_{k\in\interi^{n_{1}}\setminus\Mscr_{\omega}\atop|k|\leq rK}|\scalprod{k}{\omega}| \geq \frac{\gamma}{|k|^\tau}\ ,
\formula{frm:alphar}
$$
with  $\gamma>0$ and $\tau>n_1$.

\smallskip
\noindent
We remark that hypothesis~(i) allows us to skip the so-called {\it
geometric part} of Nekhoroshev's theorem.  The general case can be
recovered using a suitable adaptation of the geometric part of the
Nekhoroshev's theorem (see, e.g.,~\dbiref{Giorgilli-2003}).
Furthermore, assumption~(iii) implies that the Fourier expansion of
$f_0$ must have the special form, $f_0
= \sum_{k\in\Mscr_\omega}f_{0,k}(p,z,\zconj)\,
e^{\imunit \scalprod{k}{q}}$, i.e., it contains only resonant modes.

Concerning the small parameter $\epsilon$, let us remark that the
Hamiltonian~\frmref{frm:H0} is not already in the form of a function
expanded in power series of a small parameter, e.g., in $\epsilon$.
Actually, we perform an expansion in power series exploiting the
exponential decay of the Fourier coefficients, as detailed in
subsection~\sbsref{sbs:splitting}, getting
$$
\epsilon H_1 = \sum_{j\geq1} h_j(p,q,z,\zconj)\ .
$$
We stress that in the expansion above the order in the small parameter
is encoded in the index $s$ of the function $h_s$.

The quantitative estimates concerning the canonical transformation
defining the relegation are collected in the following
\proposition{pro:fondamentale}
Let the Hamiltonian~\frmref{frm:H0} satisfy the
hypotheses~(i),~(ii),~(iii) and~ (iv) above.  Take three positive
integers $K\geq1$, $r\geq1$ and $L\geq1$ and assume
$$
\frac{9 r^2 L\Xi
\mu G }{\alpha_r} \leq \frac{1}{2^{7}}\ ,
\quad\hbox{and}\quad
\eta=\frac{\epsilon r^4 A}{\alpha^2_r}+4 e^{-K \sigma/2}\leq\frac{1}{2}\ ,
\formula{frm:eqeta}
$$
with $\alpha_r$ as in $(5)$,
$$
A= 2^{21} \,\Xi^2 \left( \frac{1+e^{-\sigma/2}}{1-e^{-\sigma/2}} \right)^{n_1}
|H_1(p,q,z,\zconj)|_{\rho,2\sigma,R}\ ,
\quad
\Xi = \left(\frac{2}{e\rho\sigma}+
\frac{1}{R^{2}}\right)\ .
\formula{frm:eqA}
$$
Then:

\item{(i)}
there exist a truncated sequence of generating functions
$\Xscr^{(r)}=\{\Xscr_1,\,\ldots,\,\Xscr_r,\,0,\,\ldots\}$ that
transforms the Hamiltonian into
$$
H^{(r)} = h_0 + \mu f_0 + Z_1+\ldots+Z_r+\Rscr^{(r+1)}\ ,
$$
where the functions $Z_1,\,\ldots,\,Z_r$ are in {\it normal form},
i.e., they commute with $h_0$.  The term $\Rscr^{(r+1)}$ is the
reminder of the transformation, i.e., the collection of all the terms
that are at least of order $r+1$;

\item{(ii)}
the generating sequence defines an analytic canonical transformation on
the domain $\Dscr_{\frac{3}{4}(\rho,\sigma,R)}$ such that
$$
\Dscr_{\frac{5}{8}(\rho,\sigma,R)} \subseteq \telchi_{\Xscr}\Dscr_{\frac{3}{4}(\rho,\sigma,R)} \subseteq \Dscr_{\frac{7}{8}(\rho,\sigma,R)}\ ,\quad
\Dscr_{\frac{5}{8}(\rho,\sigma,R)} \subseteq\telchi_{\Xscr}^{-1}\Dscr_{\frac{3}{4}(\rho,\sigma,R)} \subseteq \Dscr_{\frac{7}{8}(\rho,\sigma,R)}\ ,
$$
and moreover in $\Gscr_{\frac{3}{4}\rho}$ one has
$$
|\telchi_{\Xscr}p-p| \leq \frac{\rho}{16}\ ,\quad
|\telchi_{\Xscr}^{-1}p-p| \leq \frac{\rho}{16}\ ;
\formula{frm:stima-p}
$$

\item{(iii)} the remainder is estimated by
$$
\|\Rscr^{(r+1)}\|_{\frac{3}{4}(\rho,\sigma,R)} \leq \epsilon
\left(\frac{A}{2^{18} \Xi^2}\right)
\eta^r
\ .
$$
\endclaim

\noindent
Explicit bounds for $\mu$ and $\epsilon$ are easily obtained
from~\frmref{frm:eqeta} using also~\frmref{frm:alphar}.

\remark
The remainder $\Rscr^{(r+1)}$ is composed of two kind of terms,
the ones corresponding to the truncation at a finite order of the
transformation itself, like in the usual normalisation process, and
the ones obtained at each order as a consequence of solving the
homological equation in a recursive manner.  This is the key
ingredient of the relegation algorithm and the main advantage with
respect to the classical Birkhoff normal form.  However, the estimate
at point~(iii) in proposition~\proref{pro:fondamentale} somehow hides
this distinction, being essentially the same as in the classical
normalisation algorithm (see, e.g.,~\dbiref{Giorgilli-2003}).  We
stress that the contribution due to the recursive solving of the
homological equation is controlled by~\frmref{frm:eqeta}.  Indeed the
smallness condition on $\mu$ allows to perform $L$ steps in the
recursive solving of the homological equations, getting a {\it small
remainder} term of order $\Oscr(\mu^{L+1})$.
\smallskip

We now should make a clear distinction between the {\it resonant} and
{\it non-resonant} systems.  In fact the formal scheme and the
rigorous estimates in proposition~\proref{pro:fondamentale} apply to
both cases, but the result concerning the effective stability time
exhibits a substantial difference.

In the resonant case, we have $n_1-\dim\Mscr_{\omega}$ approximate
first integrals of the form $\Phi = \telchi_{\Xscr^{(r)}}\Phi_0$, with
$\Phi_0 = \scalprod{\lambda}{p}$ where $\lambda\in \reali^{n_1-\dim\Mscr_{\omega}}$ and
satisfies $\lambda\,\bot\,\Mscr_\omega$.  We may assume that $|\lambda|<1$
without loss of generality.  We denote by $\Pi_{\Mscr_{\omega}}(p)$
the plane through $p$ generated by the resonance module
$\Mscr_\omega$, namely
$$
\Pi_{\Mscr_{\omega}}(p) = \bigl\{
p'\in\reali^{n_1} : p'-p\in {\rm span}(\Mscr_\omega)
\bigr\}\ .
$$
Using the language introduced by
Nekhoroshev\bibref{Nekhoroshev-1977}\bibref{Nekhoroshev-1979}, we call
plane of {\it fast drift} the plane $\Pi_{\Mscr_{\omega}}(p(0)) = p(0)
+ {\rm span}(\Mscr_\omega)$ and {\it deformation} the non-linear
contributions to the functions $\Phi$, which cause the orbit to
oscillate around the plane of fast drift. Finally, we say that the
remainder generates a {\it noise} that may cause a slow motion of the
orbit in a direction transversal to the plane of fast drift; we call
this slow motion {\it diffusion}.  In contrast with Birkhoff normal
form, in the relegation procedure we do not have a strict control on
the term $f_0$ and the dynamics described by $(z,\zconj)$.  Therefore,
in order to obtain a result concerning the long time stability of the
approximate first integrals we have to assume that the evolution of
$(z,\zconj)$ is confined.

Proposition~\proref{pro:fondamentale} allows to obtain a bound for the
diffusion time but, without further assumptions, we do not have any
control on the fast drift along the resonant plane. Thus, in the
resonant case, we can get an estimate of the stability time that must
be combined with some {\it a priori} bound confining the orbits.  The
formal statement of the local stability is given by the following

\lemma{lem:locale}
With the same hypotheses of proposition~\proref{pro:fondamentale}, the
following statement holds true: if $(p(t),q(t),z(t),\zconj(t))$ is an
orbit lying in the domain $\Dscr$ for $t\in[\tau^{-},\,\tau^{+}]\subset\reali$,
with $\tau^{-} < 0 < \tau^{+}$, then one has
$$
\dist\Bigl(p(t),\Pi_{\Mscr_{\omega}}(p(0))\Bigr) < \frac{\rho}{2}
$$
for all $t\in[\tau^{-},\,\tau^{+}]\cap[-t^*,\,t^*]$, with
$$
t^* = \frac{2^{12} e \rho\sigma \Xi^2}{A \epsilon \eta^r}\ ,
$$
where $\eta$, $A$ and $\Xi$ are defined as in~\frmref{frm:eqeta} and~\frmref{frm:eqA},
respectively.
\endclaim

\noindent
Let us stress that an extension of the present result can be achieved
constructing a suitable covering of boxed domains and studying the
so-called geography of resonances.  This constitutes the geometric
part of the Nekhoroshev's theorem and we refer
to~\dbiref{Giorgilli-2003} for a detailed exposition on the subject.
Such an extension requires a broad discussion, but does not contain any
essential modification with respect to the Nekhoroshev's theorem.
Moreover, as the aim of this paper is to give a rigorous support for
the relegation algorithm, we decided to state the results in the
simplest framework.

In the non-resonant case proposition~\proref{pro:fondamentale} is
enough to guarantee an exponentially long-time stability for a
suitable open set of initial data of the actions $p$.  The results for
non-resonant systems about the effective stability time is given by
the following

\theoremtx{teo:stabilita}{(non-resonant)}
Let the Hamiltonian~\frmref{frm:H0} satisfy the same hypotheses of
proposition~\proref{pro:fondamentale} and assume that the frequency
vector $\omega$ is non-resonant, i.e., $\Mscr_\omega=\{0\}$.  There
exist positive real constants $\epsilon^*$, $\mu^{*}$ and $\Tscr$ such that
the following statement holds true: if $\epsilon<\epsilon^*$, for
every orbit $(p(t),q(t),z(t),\zconj(t))$ lying in the domain $\Dscr$
at $t=0$ with $(z(t),\zconj(t))\in\Bscr$ for $t\in[\tau^{-},\,\tau^{+}]\subset\reali$,
with $\tau^{-} < 0 < \tau^{+}$ one has
$$
\dist\Bigl(p(t), p(0)\Bigr) < \frac{\rho}{2}
$$
for all $t\in[\tau^{-},\,\tau^{+}]\cap[-t^*,\,t^*]$, with
$$
t^* \leq \frac{\Tscr}{\epsilon}
\exp\left(\biggl(\frac{\gamma^2}{2\epsilon e A K^{2\tau}}\biggr)^{\frac{1}{4+2\tau}}\right)\ ,
$$
where $A$ is defined as in~\frmref{frm:eqA}.
\endclaim

\noindent
Explicit estimates for the values of $\mu^*$, $\epsilon^*$ and $\Tscr$
can be found in the proof.

Let us remark that if one is interested in actual applications to
physical models, in general, the purely analytic estimates turn out to
be too pessimistic and actually unpractical.  Nevertheless, the use of
Lie transforms provides a constructive normalisation algorithm that
can be easily translated into a recursive scheme of estimates.  Thus,
using computer algebra in order to perform high-order perturbation
expansions, one can produce estimates on the long-time stability for
{\it realistic} models.  For instance, see,
e.g.,~\dbiref{SanLhoLem-2014} and \dbiref{SanLhoLem-2015} for the
study of the effective resonant stability in the spin-orbit problem
and~\dbiref{GioLocSan-2009}, \dbiref{SanLocGio-2011}, \dbiref{SanLocGio-2013}
for the problem of the long-time stability of some of the giant
planets of the Solar system.

The paper is organised as follows. In section~\secref{sec:formale} we
reformulate the relegation algorithm via Lie transform in a suitable
way in order to translate it into a scheme of recursive estimates.
The analytic tools are reported in section~\secref{sec:tools} while
the quantitative estimate are gathered in section~\secref{sec:stime}
where we also report the proof of
proposition~\proref{pro:fondamentale}, lemma~\lemref{lem:locale} and
theorem~\thrref{teo:stabilita}.

\section{sec:formale}{The relegation algorithm}
The basis of our construction is the well known Birkhoff normal form
for a Hamiltonian system.  We follow a quite standard approach, see,
e.g.,~\dbiref{Giorgilli-2003} for a detailed discussion
or~\dbiref{GioLocSan-2014} for an exposition in a quite similar
framework.  The main problem we have to face is the presence of the
so-called small divisors, thus we need to split the perturbation in
such a way that, at every step of the normalisation process, we take
into account only a finite number of Fourier harmonics.

In the following, a special role will be played by those functions
which have a finite Fourier representation.  Thus let us introduce
some particular classes of functions.
\definition{def:classifunzioni}
Given $K_1,K_2\in \naturali$, an analytic function $f(p,q,z,\zconj)$ is said to be of class $\Pscr_{K_1,K_2}$ if 
$$
f=\sum_{k\leq K_2} c_k(p,z,\zconj) \, e^{\imunit \scalprod{k}{q}}\ ,
$$
and $c_k\in\complessi$ such that $c_k\neq0$ if and only if $k=k'+k''$, with $k'\in\Mscr_\omega$ and
$|k''|\leq K_1\,$.
\endclaim

Let us stress that the classes of functions introduced above will be
useful for both the formal scheme and the quantitative estimates.
Indeed $K_1$ and $K_2$ control the accumulation of the small divisors
and trigonometric degree, respectively.

\subsection{sbs:splitting}{Splitting of the Hamiltonian}
We now split the perturbation in a suitable form.  First let us
consider the Fourier expansion of the term $\epsilon H_1$
in~\frmref{frm:H0}
$$
\epsilon H_1 = \sum_{k\in\interi^{n_1}} c_k(p,z,\zconj)\, e^{\imunit \scalprod{k}{q}}\ .
$$
Then, let us pick an arbitrary positive integer $K$ and write $H_1$ in
the form
$$
\epsilon H_1 = \sum_{j\geq1} h_j(p,q,z,\zconj)\ ,
$$
where
$$
\eqalign{
h_{1}(p,q,z,\zconj) &= \sum_{0\leq |k| < K} c_k(p,z,\zconj)\, e^{\imunit \scalprod{k}{q}}\ ,\cr
\vdots\cr
h_{s}(p,q,z,\zconj) &= \sum_{(s-1)K\leq |k| < sK} c_k(p,z,\zconj)\, e^{\imunit \scalprod{k}{q}}\ .\cr
}
\formula{frm:sviluppohs}
$$
This procedure breaks down the classical scheme of series expansions
in the perturbative parameter $\epsilon$, introducing an arbitrary
quantity $K$.  However, at the end of the proof of
theorem~\thrref{teo:stabilita}, we will see that there is a natural
choice\footnote{2}{As reported
in~\dbiref{Giorgilli-2003} (pag.~86, footnote~2), choosing the
parameter $K$ by asking $e^{-K \sigma}\sim\epsilon$ is not really
convenient.  As it will be evident from the optimisation of the
parameters in the proof of theorem~\thrref{teo:stabilita}, the best
choice is $K\sim 1/\sigma$.} for $K$.

The splitting introduced above is based on the exponential decay of
the Fourier coefficients of an analytic function, as stated by the
following

\lemma{lem.decadimento}
Let $H_1(p,q,z,\zconj)$ be analytic in $\Dscr_{\rho,2\sigma,R}$ and
$$
|H_1(p,q,z,\zconj)|_{\rho,2\sigma,R}=\sup_{{\scriptstyle{p\in\Gscr_{\rho}\,,\ q\in\toro^{n_1}_{2\sigma}}}
      \atop{{\scriptstyle{(z,\zconj)\in \Bscr_{R}}}}}
\big|H_1(p,q,z,\zconj)\big|
<\infty\ .
$$
Then
$$
\|h_s\|_{\rho,\sigma,R} \leq \zeta^{s-1} F\ ,\quad s\geq1\ ,
$$
with
$$
\zeta = e^{-K \sigma/2}\ ,\quad
F=\epsilon\left( \frac{1+e^{-\sigma/2}}{1-e^{-\sigma/2}} \right)^{n_1}
|H_1(p,q,z,\zconj)|_{\rho,2\sigma,R}\ .
$$
\endclaim

\noindent
The proof of lemma~\lemref{lem.decadimento} is straightforward, see,
e.g.,~lemma~5.2 in~\dbiref{Giorgilli-2003}.

\subsection{sbs:formalscheme}{Formal scheme}
Let us pick two positive integers $K$, $K'$ and write the
Hamiltonian~\frmref{frm:H0} as
$$
H^{(0)} = h_0 + \mu f_0 + h_1 + h_2 +
\ldots +  h_s+\ldots\ ,
\formula{frm:espansioneH0}
$$
with $h_0=\scalprod{\omega}{p}$, $f_0 \in \Pscr_{0,K'}$ and
$h_s \in \Pscr_{sK,sK}$.  We remark that with this splitting of the
perturbation, the terms $h_s$ are of order $s$ in some small
parameter, precisely they are of order $\Oscr(\epsilon\, \zeta^{s-1})$.

We look for a sequence of generating functions
$\Xscr^{(r)}=\{\Xscr_1,\,\ldots,\,\Xscr_r,\,0,\,\ldots\}$, with $r$ arbitrary
positive integer, and a function $Z^{(r)}=Z_0+\ldots+Z_r$ such that
$$
\telchi_{\Xscr^{(r)}} Z^{(r)} = h_0 + \mu f_0 + h_1 + \ldots + h_r + \Qscr^{(r+1)}\ .
\formula{frm:formanormale-r}
$$
The functions $Z_0,\,\ldots,\,Z_r$ must be determined so as to be in
{\it normal form}, namely, they must commute with $h_0$.  Instead, the
term $\Qscr^{(r+1)}$ is the unrelegated remainder, namely
$$
\telchi_{\Xscr^{(r)}} Z^{(r)} - \sum_{s=0}^{r} E_s Z^{(r)} = \Qscr^{(r+1)}\ .
$$
The remainder, being a term of order $r+1$ in the small parameter, can
be considered as a {\it small} term. Moreover, the asymptotic
character of the transformation allows to progressively lower the
influence of the remainder. However, it is important to recall that
the resulting series are actually divergent, thus the sequence of
canonical transformations must be {\it finite}.

Again, we emphasise that the relegation algorithm represents a {\it
variazione} of the classical Birkhoff normal form.  The difference is
the special role played by the term $f_0$.  Thus, we introduce the
additional parameter $L$ so as to take into account the peculiar
character of $f_0$.  The role of $L$ will be clear from the definition
of the generating functions.

For $r=1$, we have to solve the following equations
$$
\eqalign{
Z_{1,0} - \lie{h_0} \Xscr_{1,0} &= h_1\ ,\cr
Z_{1,j} - \lie{h_0} \Xscr_{1,j} &= \lie{\mu f_0} \Xscr_{1,j-1}\ ,\quad j=1,\,\ldots,\,L\ .\cr
}
$$
We define the generating function $\Xscr_1$ as
$$
\Xscr_1 = \Xscr_{1,0}+\Xscr_{1,1}+\ldots+\Xscr_{1,L}\ .
$$
Explicit expressions for $\Xscr_1$ are easy to obtain.  We  report
here the sole expression of $\Xscr_{1,0}$, that follows
from~\frmref{frm:sviluppohs}
$$
\Xscr_{1,0}(p,q,z,\zconj) = \sum_{k\in\interi^{n_1}\setminus \Mscr_\omega \atop 0\leq |k| < K} \frac{c_k(p,z,\zconj)}{\imunit\scalprod{k}{\omega}}
\, e^{\imunit \scalprod{k}{q}}\ ,
$$
similar expressions can be obtained for $\Xscr_{1,j}$ with $j=1,\ldots L$.
Let us remark that
$$
\Xscr_{1,j} \in \Pscr_{K,K+jK'}\ ,\quad \Xscr_{1} \in \Pscr_{K,K+LK'}\ .
$$
The relegated term $Z_1$ is given by
$$
Z_1 = Z_{1,0}+Z_{1,1}+\ldots+Z_{1,L}\ .
$$
Thus for $r=1$ we get
$$
Z_{1} - \lie{Z_0} \Xscr_{1} = h_1 + \lie{\Xscr_{1,L}} \mu f_0\ .
$$
Proceeding by induction, for $r>1$, for all orders $s \leq r$ we need to solve
$$
\eqalign{
Z_{s,0} - \lie{h_0} \Xscr_{s,0} &= \Psi_r\ ,\cr
Z_{s,j} - \lie{h_0} \Xscr_{s,j} &= \lie{\mu f_0} \Xscr_{s,j-1}\ ,\quad j=1,\,\ldots,\,L\ ,\cr
}
$$
where
$$
\eqalign{
\Psi_{1} &= h_1\ ,\cr
\Psi_{2} &= h_2 - \frac{1}{2} \left( \lie{\Xscr_1}h_{1} + E_{1}Z_1\right)- \frac{1}{2} \lie{\Xscr_1}
\lie{\Xscr_{1,L}} \mu f_0
-\lie{\Xscr_{1,L}}\mu f_0
\ ,\cr
\Psi_{s} &= h_s - \sum_{j=1}^{s-1} \frac{j}{s} \left( \lie{\Xscr_j}h_{s-j} + E_{s-j}Z_j\right) \cr
&\quad-\sum_{j=1}^{s-1} \frac{j}{s} \lie{\Xscr_j}
\left( \lie{\Xscr_{s-j,L}} \mu f_0 - \lie{\Xscr_{s-j-1,L}} \mu f_0 \right)
-\lie{\Xscr_{s-1,L}}\mu f_0 \ ,\quad s>2\ ,\cr
}
$$
with $\Psi_{s} \in \Pscr_{sK,s(K+LK')}\,$.  

Let us remark that we
cannot let either $L$ or $r$ go to infinity, as the radius of
convergence shrinks to zero; this is in agreement with the classical
Birkhoff normal form.

\section{sec:tools}{Analytic tools}
In order to make the paper self-contained, we report in this section
all the technical tools needed to prove
proposition~\proref{pro:fondamentale}, lemma~\lemref{lem:locale} and
theorem~\thrref{teo:stabilita}.  Again we recall that we use the
shorthand notations $|\cdot|_{\alpha}\,$ and $\|\cdot\|_{\alpha}\,$
for $|\cdot|_{\alpha(\rho,R)}$ and
$\|\cdot\|_{\alpha(\rho,\sigma,R)}\,$.

\subsection{app:stima-serie-di-Lie}{Estimates for multiple Poisson brackets}
Some Cauchy estimates on the derivatives in the restricted domains
will be useful.

\lemma{lem:stima-derivata-Lie}
Let $d\in\reali$ such that $0<d<1$ and $g$ be an analytic function with bounded norm
$\|g\|_{1}$.  Then one has
$$
\left\|\parder{g}{p_j}\right\|_{1-d}\leq
\frac{\left\|g\right\|_{1}}{d\rho}\ ,
\qquad
\left\|\parder{g}{q_j}\right\|_{1-d}\leq
\frac{\left\|g\right\|_{1}}{e d \sigma}\ ,
\qquad
\left\|\parder{g}{z_j}\right\|_{1-d}\leq
\frac{\left\|g\right\|_{1}}{d R}\ .
\formula{frm:stime-Cauchy}
$$
Of course, the latter inequality holds true also by replacing $z_j$
with $\imunit\bar z_j\,$.
\endclaim
\noindent
The proof of lemma~\lemref{lem:stima-derivata-Lie} is straightforward and it
is left to the reader.

\lemma{lem:stima-parentesi-Poisson}
Let $d,d'\in\reali$ such that $d>0$, $d'\geq0$ and $d+d^{\prime}<1$, and
$g,\,g^{\prime}$ be two analytic functions with bounded norms
$\|g\|_{1-d-d^{\prime}}$ and $\|g^{\prime}\|_{1-d^{\prime}}$, respectively. Then,
for all $\delta\in\reali$ such that $\delta>0$ and $d+d^{\prime}+\delta<1$ one
has
$$
\left\|\{g,g^{\prime}\}\right\|_{1-d-d^{\prime}-\delta}\leq
\frac{\Xi}{(d+\delta)\delta}
\left\|g\right\|_{1-d-d^{\prime}}\left\|g^{\prime}\right\|_{1-d^{\prime}}\ ,
\quad 
\Xi = \left(\frac{2}{e\rho\sigma}+
\frac{1}{R^{2}}\right)\ .
\formula{frm:stimapois}
$$
\endclaim
\proof
We separately consider the parts of the Poisson bracket involving the
$(p,q)$ variables and the $(z,\zconj)$ ones.  For the first
part we have
$$
\eqalign{
&\Bigg\|\sum_{j=1}^{n_1}\left(
\frac{\partial g}{\partial q_j}\frac{\partial g^{\prime}}{\partial p_j}-
\frac{\partial g}{\partial p_j}\frac{\partial g^{\prime}}{\partial q_j}
\right)\Bigg\|_{1-d-d^{\prime}-\delta}
\leq
\cr
&\qquad\qquad\qquad
\sum_{k\in\interi^{n_1}}\sum_{k^{\prime}\in\interi^{n_1}}\Bigg[\Bigg(
\frac{|k|\left|g_{k}\right|_{1-d-d^{\prime}}
\left|g^{\prime}_{k^{\prime}}\right|_{1-d^{\prime}}}{(d+\delta)\rho}
\cr
&\qquad\qquad\qquad
\phantom{\sum_{k\in\interi^{n_1}}\sum_{k^{\prime}\in\interi^{n_1}}\Bigg[\Bigg(}
+\frac{\left|g_{k}\right|_{1-d-d^{\prime}}
|k^{\prime}|\left|g^{\prime}_{k^{\prime}}\right|_{1-d^{\prime}}}{\delta\rho}
\Bigg)\ee^{(|k|+|k^{\prime}|)(1-d-d^{\prime}-\delta)\sigma}\Bigg]
\cr
&\qquad\qquad\qquad
\le\frac{2}{e\rho\sigma}\frac{1}{(d+\delta)\delta}
\left\|g\right\|_{1-d-d^{\prime}}\left\|g^{\prime}\right\|_{1-d^{\prime}}\ ,
}
\formula{frm:stimavel}
$$
being $g_{k^{\prime}}^{\prime}=g_{k^{\prime}}^{\prime}(p,z,\zconj)$ as in the expansion~\frmref{frm:funz}.  Here the Cauchy
estimate~\frmref{frm:stime-Cauchy} and the elementary inequality
$a\ee^{-ab}\le 1/(\ee b)\,$, for positive $a$ and $b$ have been used.

Let us now focus on the second part of the Poisson bracket.  For every
point $(p,z,\zconj)\in\Gscr_{(1-d-d^{\prime}-\delta)\rho}
\times \Bscr_{(1-d-d^{\prime}-\delta)R}$ and for all pairs of vectors
$k,\,k^{\prime}\in\interi^{n_1}$, we introduce an auxiliary function
$$
W_{(p,z,\zconj);k,k^{\prime}}(t) =
g_k\Bigl(p,\,z-t\parder{g_{k^{\prime}}^{\prime}}{(\zconj)},
\,\zconj+t\parder{g_{k^{\prime}}^{\prime}}{z}\Bigr)\ .
$$
Since $g_k$ is analytic on $\Gscr_{(1-d-d^{\prime})\rho}\times
\Bscr_{(1-d-d^{\prime})R}\,$, then
$W_{(p,z,\zconj);k,k^{\prime}}$ is analytic for
$|t|\le \bar t\,$, with
$$
{\bar t}=\frac{\delta R}{{\displaystyle \max_{1\le j\le n_2}}
\left\{\left|\parder{g_{k^{\prime}}^{\prime}}{z_j}\right|_{1-d-d^{\prime}-\delta}\,,\,
\left|\parder{g_{k^{\prime}}^{\prime}}{(\zconj_j)}
\right|_{1-d-d^{\prime}-\delta}\right\}}\ .
$$
Thus, by the Cauchy's estimate we get
$$
\big|\left\{g_{k},g_{k^{\prime}}^{\prime}\right\}\big|_{1-d-d^{\prime}-\delta}\le
\left| \frac{\der}{\der t} W_{(p,z,\zconj);k,k^{\prime}}(t)\Bigm|_{t=0}
\right|_{1-d-d^{\prime}-\delta}
\le
\frac{\left|g_{k}\right|_{1-d-d^{\prime}}}{\bar t}\ .
$$
By the definition~\frmref{frm:normafou} of the norm, we get
$$
\Bigg\|\sum_{j=1}^{n_2}\left(
\frac{\partial g}{\partial (\zconj_j)}
\frac{\partial g^{\prime}}{\partial z_j}-
\frac{\partial g}{\partial z_j}
\frac{\partial g^{\prime}}{\partial (\zconj_j)}
\right)\Bigg\|_{1-d-d^{\prime}-\delta}\le
\frac{1}{R^2}
\frac{\left\|g\right\|_{1-d-d^{\prime}}\left\|g^{\prime}\right\|_{1-d^{\prime}}}
{(d+\delta)\delta}\ .
\formula{frm:stimasec}
$$
The wanted inequality~\frmref{frm:stimapois} follows by adding
up~\frmref{frm:stimavel} and~\frmref{frm:stimasec}.
\endproof

\lemma{lem:stima-derivata-Lie-j}{
Let $d,d'\in\reali$ such that $d>0$, $d'\geq0$ and $d+d^{\prime}<1$ and $\Xscr,\,g$ be
two analytic functions with bounded norms $\|\Xscr\|_{1-d^{\prime}}$
and $\|g\|_{1-d^{\prime}}\,$, respectively. Then, for $j\ge 1$, we
have
$$
\left\|\Lie^{j}_{\Xscr}g\right\|_{1-d-d^{\prime}}
\le\frac{j!}{\ee^{2}}
\left(e^2\frac{\Xi}{d^2}\right)^{j}
\|\Xscr\|^{j}_{1-d^{\prime}}\|g\|_{1-d^{\prime}}\ ,
\quad
\Xi = \left(\frac{2}{e\rho\sigma}+
\frac{1}{R^{2}}\right)\ .
$$
}\endclaim

\proof
For $j\ge 1$ let $\delta=d/j$.  By repeated application of
lemma~\lemref{lem:stima-parentesi-Poisson}, we get the recursive chain
of inequalities
$$
\eqalign{
\left\|\Lie^{j}_{\Xscr}g\right\|_{1-d-d^{\prime}}
&\leq
\frac{\Xi}{j\delta^{2}}\|\Xscr\|_{1-d^{\prime}}
\left\|\Lie^{j-1}_{\Xscr}g\right\|_{1-d^{\prime}-(j-1)\delta}
\cr
&\leq\ldots
\cr
&\leq\frac{j!}{\ee^{2}}
\left(e^2\frac{\Xi}{d^2}\right)^{j}
\|\Xscr\|^{j}_{1-d^{\prime}}\|g\|_{1-d^{\prime}}\ .
\cr
}
$$
In the last row, we used the trivial inequality $j^{j}\leq
j!\,\ee^{j-1}$, holding true for $j\ge 1\,$.
\endproof

\subsection{analiticita}{Analyticity of Lie transform}
We report here two main results concerning the analyticity of the Lie
transform.

\proposition{prop:1}
Let the generating sequence $\Xscr=\{\Xscr_s\}_{s\geq1}$ be analytic on the
domain $\Dscr_{\rho,\sigma,R}\,$, and assume
$$
\|\Xscr_s\|_{\rho,\sigma,R} \leq \frac{b^{s-1}}{s}G\ ,
\formula{frm:stima-generatrice}
$$
with $b,G\in\reali$ such that $b\geq0$ and $G>0$. Then, for every positive $d<1/2$ the
following statement holds true: if the condition
$$
e^2 \frac{\Xi}{d^2}G+b\leq
\frac{1}{2}\ ,
\quad
\Xi = \left(\frac{2}{e\rho\sigma}+
\frac{1}{R^{2}}\right)
\formula{frm:condizione-Xi}
$$
is satisfied, then the operator $\telchi_{\Xscr}$ and its inverse
$\telchi_{\Xscr}^{-1}$ define an analytic canonical transformation on
the domain $\Dscr_{(1-d)(\rho,\sigma,R)}$ with the properties
$$
\vcenter{\openup1\jot\halign{%
${#}$
&$\subset\ ${#}
&${#}$
&$\subset\ ${#}
&${#}$
\cr
\Dscr_{(1-2d)(\rho,\sigma,R)} &&\telchi_{\Xscr}\Dscr_{(1-d)(\rho,\sigma,R)} &&\Dscr_{\rho,\sigma,R}\ ,\cr
\Dscr_{(1-2d)(\rho,\sigma,R)} &&\telchi_{\Xscr}^{-1}\Dscr_{(1-d)(\rho,\sigma,R)} &&\Dscr_{\rho,\sigma,R}\ .\cr
}}
$$
\endclaim

The proof of the proposition is based on the following
\lemma{lemma:1}
Let a function $f$ and the generating sequence $\Xscr=\{\Xscr_s\}_{s\geq1}$ be
analytic on the domain $\Dscr_{\rho,\sigma,R}\,$, and assume that
$\|f\|_{\rho,\sigma,R}$ is finite.  Let the generating sequence
satisfy~\frmref{frm:stima-generatrice} and assume
that~\frmref{frm:condizione-Xi} holds true.  Then the series
$\telchi_{\Xscr}f$, $\telchi_{\Xscr}^{-1}f$, $\telchi_{\Xscr}p$,
$\telchi_{\Xscr}^{-1}p$, $\telchi_{\Xscr}q$ and
$\telchi_{\Xscr}^{-1}q$ are absolutely convergent on
$\Dscr_{(1-d)(\rho,\sigma,R)}\,$, and for any integer $r>0$ one has

\item{(i)} the operators $E_s$ and $D_s$ are estimated by
$$
\eqalign{
\|E_s f\|_{(1-d)(\rho,\sigma,R)} &\leq \left(
e^2 \frac{\Xi}{d^2}G+b
\right)^{s-1}
 \frac{\Xi}{d^2}G
\|f\|_{\rho,\sigma,R}\ ,\cr
\|D_s f\|_{(1-d)(\rho,\sigma,R)} &\leq \left(
e^2 \frac{\Xi}{d^2}G+b
\right)^{s-1}
\frac{\Xi}{d^2}G
\|f\|_{\rho,\sigma,R}\ ;\cr
}
$$

\item{(ii)} the Lie transform and its inverse are estimated by
$$
\|\telchi_{\Xscr}f\|_{(1-d)(\rho,\sigma,R)} \leq 2 \|f\|_{\rho,\sigma,R}\ ,
\quad
\|\telchi_{\Xscr}^{-1}f\|_{(1-d)(\rho,\sigma,R)} \leq 2 \|f\|_{\rho,\sigma,R}\ ;
$$

\item{(iii)} the remainder of a $r$-th order truncated transformation is estimated by
$$
\eqalign{
\|\telchi_{\Xscr} f - \sum_{s=0}^{r}E_s f\|_{(1-d)(\rho,\sigma,R)} &\leq
\frac{2}{e^2}\left(
e^2\frac{\Xi}{d^2}G+b
\right)^{r+1}
\|f\|_{\rho,\sigma,R}\ ,\cr
\|\telchi_{\Xscr}^{-1} f - \sum_{s=0}^{r}D_s f\|_{(1-d)(\rho,\sigma,R)} &\leq
\frac{2}{e^2}\left(
e^2\frac{\Xi}{d^2}G+b
\right)^{r+1}
\|f\|_{\rho,\sigma,R}\ ;\cr
}
$$

\item{(iv)} the change of coordinates is estimated by
$$
\vcenter{\openup1\jot\halign{%
${#}$ & ${#}$ & \quad${#}$ & ${#}$ \cr
\|\telchi_{\Xscr}p-p\|_{(1-d)(\rho,\sigma,R)} &\leq \frac{1}{2e^2} d\rho\ ,
&\|\telchi_{\Xscr}^{-1}p-p\|_{(1-d)(\rho,\sigma,R)} &\leq \frac{1}{2e^2} d\rho\ ,\cr
\|\telchi_{\Xscr}q-q\|_{(1-d)(\rho,\sigma,R)} &\leq \frac{1}{2e} d\sigma\ ,
&\|\telchi_{\Xscr}^{-1}q-q\|_{(1-d)(\rho,\sigma,R)} &\leq \frac{1}{2e} d\sigma\ ,\cr
\|\telchi_{\Xscr}z-z\|_{(1-d)(\rho,\sigma,R)} &\leq \frac{1}{2e^2} d R\ ,
&\|\telchi_{\Xscr}^{-1}z-z\|_{(1-d)(\rho,\sigma,R)} &\leq \frac{1}{2e^2} d R\ .\cr
}}
$$
\endclaim

The proof of the previous lemma and proposition~\proref{prop:1} are
just a straightforward adaptation of proposition~4.3 and lemma~4.4
in~\dbiref{Giorgilli-2003}.

\section{sec:stime}{Quantitative estimates}
We now translate the formal scheme introduced in
section~\secref{sec:formale} into recursive estimates.

The quantitative estimates for the truncated sequence of generating
functions are collected in the following
\lemma{lem:stima-generatrici}
Consider the Hamiltonian~\frmref{frm:espansioneH0} and let
$\|h_s\|_1\leq\zeta^{s-1}F$, for some real $F>0$ and $\zeta\geq0$, and
assume $\|f_0\|_{\rho,\sigma,R}\leq G$.  Furthermore, for any positive
$d<1$, let $\mu$ satisfy the smallness condition
$$
\frac{9 r^2 L\Xi
\mu G }{\alpha_r d^2} \leq \frac{1}{2}\ ,
\quad
\Xi = \left(\frac{2}{e\rho\sigma}+
\frac{1}{R^{2}}\right)\ ,
$$
with $\alpha_r$ as in~\frmref{frm:alphar}.  Then the truncated
sequence of generating functions $\Xscr^{(r)}$ that gives the
Hamiltonian the normal form~\frmref{frm:formanormale-r} satisfies
$$
\|\Psi_s\|_{(1-d)(\rho,\sigma,R)}\leq \frac{b^{s-1}}{s} F\ ,\quad
\|\Xscr_s\|_{(1-d)(\rho,\sigma,R)}\leq \frac{b^{s-1}}{s} \frac{2F}{\alpha_r}\ ,
$$
with
$$
b = 4\left(\frac{2^7 r^4  F\, \Xi^2}{\alpha_r^2 d^4}
+ \zeta\right)\ .
$$
\endclaim

\proof
In order to produce recursive estimates of the norm of Poisson
brackets, we need to define a suitable sequence of restrictions of the
domain.  Fix the final restriction, $d$, and define the following
sequence
$$
d_i = \frac{i}{3r} d\ ,\quad
i=1,\,\ldots,\,3r\ ,
\formula{frm:restringimenti}
$$
where $r$ is the maximum relegation order.  We now look for the
estimates of $\|\Psi_s\|_{1-d_{3s-2}}$ and $\|\Xscr_s\|_{1-d_{3s-1}}$.
Again, here we omit the $(\rho,\sigma,R)$ in the subscript of the norms.

For $s=1$ we have
$$
\|\Psi_1\|_{1-d_1} \leq F\ ,
$$
from which we immediately get
$$
\|\Xscr_{1,0}\|_{1-d_1} \leq \frac{F}{\alpha_r}\ .
$$
Iterating the previous estimates, for $l\geq1$ we get
$$
\eqalign{
\|\Xscr_{1,l}\|_{1-d_1-l\frac{d_2-d_1}{L}} &\leq
\frac{1}{\alpha_r}
\frac{\Xi}{\left( d_1+l\frac{d_2-d_1}{L}\right)\frac{d_2-d_1}{L}}
\|\Xscr_{1,l-1}\|_{1-d_1-(l-1)\frac{d_2-d_1}{L}} \|\mu f_0\|_1\cr
&\leq
\frac{9 r^2 L\,\Xi}{\alpha_r d^2}\,
\|\Xscr_{1,l-1}\|_{1-d_1-(l-1)\frac{d_2-d_1}{L}} \|\mu f_0\|_1\ ,
}
$$
where we used the elementary inequality
$$
\left(d_{1}+l\frac{d_2-d_1}{L}\right) \frac{d_2-d_{1}}{L} \geq d_{1}\frac{d_2-d_{1}}{L} \geq \frac{d^2}{9 r^2 L}\ .
$$
Plugging in the condition
$$
\frac{9 r^2 L\Xi
\mu G }{\alpha_r d^2} \leq \frac{1}{2}\ ,
$$
we get
$$
\|\Xscr_1\|_{1-d_2} \leq
2\frac{F}{\alpha_r}\ .
$$
Summing up, for $s=1$, we have
$$
\|\Psi_1\|_{1-d_1} \leq F\ ,\quad
\|Z_1\|_{1-d_1} \leq F\ ,\quad
\|\Xscr_1\|_{1-d_2} \leq
2\frac{F}{\alpha_r}\ .
$$

We now look for two real sequences $\{\eta_s\}_{1\leq s\leq r}$ and
$\{\theta_{s,j}\}_{0\leq s\leq r\,,\ 1\leq j\leq r}$ such that
$$
\|\Psi_s\|_{1-d_{3s-2}}\leq \eta_s F\ ,\quad
\|E_s Z_j\|_{1-d_{3(s+j)-2}} \leq \theta_{s,j} F\ ,
$$
and we remark that we may choose $\eta_1=1$ and set
$\theta_{0,j}=\eta_j$.

For $s\geq2$ we get
$$
\eqalign{
\|\Psi_s\|_{1-d_{3s-2}} &\leq
\zeta^{s-1}F +\sum_{j=1}^{s-1} \frac{j}{s} \biggl(
\frac{\Xi}{d_{3s-2}(d_{3s-2}-d_{3j-1})} 2\eta_j \frac{F}{\alpha_r}
\zeta^{s-j-1}+\theta_{s-j,j}\biggr)F\cr
&\qquad+
\sum_{j=1}^{s-1} \frac{j}{s}
\left(\frac{\Xi}{(d_{3s-2}-d_{3j-1})(d_{3s-2}-d_{3(s-j)})} \right.\cr
&\qquad\qquad \left.\frac{\Xi}{d_{3(s-j)}(d_{3(s-j)}-d_{3(s-j)-1})}\right)
8\eta_{j}\eta_{s-j}\frac{F^2}{\alpha_r^2}\cr
&\qquad+
\frac{\Xi}{d_{3s-2}(d_{3s-2}-d_{3s-4})} 2\eta_{s-1}\frac{F^2}{\alpha_r}\ ,\cr
}
$$
and, by the definition of $E_s$, for $s\geq1$ we have
$$
\|E_s Z_j\|_{1-d_{3(s+j)-3}} \leq
\sum_{l=1}^{s}\frac{l}{s}
\frac{2\eta_l\theta_{s-l,j}\Xi}{(d_{3(s+j)-2}- d_{3l-1})(d_{3(s+j)-2}-d_{3(s+j-l)-2})}
\frac{F^2}{\alpha_r}\ .
$$

In view of \frmref{frm:restringimenti}, for $1\leq j\leq s-1$ one has
$$
\frac{1}{(d_{3s-2})(d_{3s-2}-d_{3j-1})}\leq \frac{3 r^2}{2d^2}\ ,
$$
and
$$
\frac{1}{(d_{3s-2}-d_{3j-1})(d_{3s-2}-d_{3(s-j)})\,d_{3(s-j)}(d_{3(s-j)}-d_{3(s-j)-1})}   \leq
\frac{27 r^4}{2 d^4}\ ,
$$
while for $s\geq 1$, $j\leq s$ and $1\leq l\leq s$ one gets
$$
\frac{1}{(d_{3(s+j)-2}- d_{3l-1})(d_{3(s+j)-2}-d_{3(s+j-l)-2})}\leq \frac{3 r^2}{2d^2}\ .
$$

Thus the sequences $\{\eta_s\}$ and $\{\theta_{s,j}\}$ may be defined
as
$$
\eqalign{
\eta_s &= \zeta^{s-1} + \frac{C_r}{s} \sum_{j=1}^{s-1} j \eta_j \zeta^{s-j-1}  + \frac{C_r}{s} \sum_{j=1}^{s-1} j \eta_j \eta_{s-j}
+ \frac{1}{s} \sum_{j=1}^{s-1} j \theta_{s-j,j}\ ,\cr
\theta_{s,j} &= \frac{C_r}{s} \sum_{l=1}^{s} l\eta_l\theta_{s-l,l}\ , \cr}
$$
with
$$
C_r = \frac{2^7 r^4  F\, \Xi^2}{\alpha_r^2 d^4}\ .
$$

We are led to study the behaviour of the two sequences with initial
values $\eta_1=1$ and $\theta_{0,j}=\eta_j$.  Let us remark that the
second sequence can be rewritten as $\theta_{s,j}
= \eta_j \theta_{s,1}$, hence we set $\theta_{s,1}=\theta_{s}$ and
consider the double sequence
$$
\eqalign{
\eta_s &= \zeta^{s-1} + \frac{C_r}{s} \sum_{j=1}^{s-1} j \eta_j \zeta^{s-j-1} + \frac{C_r}{s} \sum_{j=1}^{s-1} j \eta_j \eta_{s-j}
+ \frac{1}{s} \sum_{j=1}^{s-1} j \eta_j \theta_{s-j}\ ,\cr
\theta_{s} &= \frac{C_r}{s} \sum_{j=1}^{s} j\eta_j\theta_{s-j}\ , \cr}
\formula{eq:etathetas}
$$
with initial values $\eta_1=\theta_0=1$.  Subtracting from the second
equation the first one multiplied by $C_r$ we get
$$
\theta_s = C_r \eta_s - C_r \zeta^{s-1} - \frac{C_r^2}{s} \sum_{j=1}^{s-1} j \eta_j \zeta^{s-j-1}\ .
$$
Let us remark that $\theta_s$ just depends by
$\eta_1\,,\ldots\,,\eta_{s-1}$.  Replacing the previous expression in
the first of~\frmref{eq:etathetas} we get
$$
\eta_s \leq \zeta^{s-1} + C_r \sum_{j=1}^{s-1} \eta_j \eta_{s-j}\ ,
$$
and we have
$$
\eta_s \leq (C_r+\zeta)^{s-1} \nu_s\ ,
\formula{frm:etas}
$$
where the real sequence $\{\nu_s\}_{s\geq1}$ is the so-called Catalan sequence
$$
\nu_1=1\ ,\quad \nu_s =  \sum_{j=1}^{s-1} \nu_j \nu_{s-j}\leq\frac{4^{s-1}}{s}\ .
$$
The claim follows replacing $\nu_s\leq 4^{s-1}/s$
in~\frmref{frm:etas} and collecting the estimates previously
obtained.
\endproof

\prooftx{of proposition~\proref{pro:fondamentale}}
The proof is a straightforward application of the results previously
obtained, having set $d=1/8$, and it is left to the reader.
\endproof

\prooftx{of lemma~\lemref{lem:locale}}
Using the elementary
estimate
$$
|\Phi_0(t) - \Phi_0(0) | \leq
|\Phi_0(t) - \Phi(t) | +
|\Phi(t) - \Phi(0) | +
|\Phi(0) - \Phi_0(0) |\ .
$$
By means of~\frmref{frm:stima-p} we easily bound the sum of the first and
third terms, being smaller than $\rho/8$.  Coming to the second term
we have
$$
|\dot\Phi|\leq\|\telchi_{\Xscr} \poisson{\Phi_0}{\Rscr^{(r+1)}}\|_{\frac{1}{2}}
\leq
\frac{16}{e \sigma}
\left(\frac{\epsilon A}{2^{18} \Xi^2}\right)
\eta^r
\leq
\frac{\epsilon A}{2^{14} e \sigma \Xi^2} \eta^r\ .
$$
Thus
$$
|\Phi(t)-\Phi(0)|\leq
|t| \frac{A}{2^{14} e \sigma \Xi^2}\, \epsilon \eta^r\ ,
$$
which is smaller than $\rho/4$ if $|t|<t^*$ as claimed.
\endproof

The proof of theorem~\thrref{teo:stabilita} follows directly from
lemma~\lemref{lem:locale}, it is just a matter of making a clever
choice of the parameters.

\prooftx{of Theorem~\thrref{teo:stabilita}}
It remains to choose the parameters $K\geq1$ and $r\geq1$ as functions
of the parameters $\rho$, $\sigma$, $R$ and $\epsilon$ that
characterise the Hamiltonian.  The aim is to make a good choice, so
that the stability time is as large as possible.

Assuming the smallness condition $\eta<1/e$, namely
$$
\frac{\epsilon r^4 A}{\alpha^2_r}+4 e^{-K \sigma/2}\leq\frac{1}{e}\ ,
$$
that is satisfied if
$$
\frac{\epsilon r^4 A}{\alpha^2_r} \leq \frac{1}{2e}\ ,
\quad
e^{-K \sigma/2}\leq\frac{1}{8e}\ ,
$$
where $A$ is defined as in~\frmref{frm:eqA}.  Using the Diophantine
condition
$$
\alpha_r \leq \frac{\gamma}{(rK)^\tau}\ ,
$$
it is natural to choose
$$
r = \biggl(\frac{\gamma^2}{2\epsilon e A K^{2\tau}}\biggr)^{\frac{1}{4+2\tau}}\ ,
\quad
K = \left\lceil \frac{2(1+3\log 2)}{\sigma} \right\rceil\ .
$$
We have $K\ge1$ by definition, while the condition $r\geq1$ is satisfied provided
$$
\epsilon \leq \frac{\gamma^2}{2 e A K^{2\tau}}\ .
$$
The claim follows by setting $\Tscr={2^{12} e \rho \sigma \Xi^2}/{A}$.
\endproof

\acknowledgements{%
We warmly thank A. Giorgilli for helpful discussions and useful
comments. The work of~M.~S. have been partially supported by the
research program ``Teorie geometriche e analitiche dei sistemi
Hamiltoniani in dimensioni finite e infinite'', PRIN 2010JJ4KPA 009,
financed by MIUR.}

\references

\bye